\documentclass[12pt]{amsart}

\def\checkbox{\leavevmode\vbox to 9pt{\hrule \vss
	\hbox to 9pt{\vrule height 9pt \hfil\vrule height 9pt}\vss
	\hrule}\ }

\newcommand{\N}{{\mathbb N}}
\newcommand{\R}{{\mathbb R}}
\newcommand{\Z}{{\mathbb Z}}
\newcommand{\C}{{\mathbb C}}

\renewcommand{\epsilon}{\varepsilon}
\renewcommand{\phi}{\varphi}

\newtheorem{Lemma}{Lemma}
\newtheorem{Theorem}[Lemma]{Theorem}

\newtheorem{Definition}[Lemma]{Definition}

\newtheorem{Conjecture}{Conjecture}

\begin{document}

\address{Department of Mathematics, Pennsylvania State University, University Park, PA $16802$, USA.}
\author{Alexander Borisov}
\title[Singularities, Factorials and RH]{Quotient singularities, integer ratios of factorials and the Riemann Hypothesis}

\begin{abstract} The goal of this paper is to reveal a close connection between the following three subjects that have not been studied together in the past:

\begin{enumerate}
\item terminal and canonical cyclic quotient singularities;

\item integer ratios of factorials;

\item Nyman's approach to the Riemann Hypothesis.

\end{enumerate}

In particular, we notice that the constructions of P.A. Picon are relevant for the study of singularities and possibly the Riemann Hypothesis. The list of the 29 stable quintuples of Mori-Morrison-Morrison coincides, up to the choice of notation, with the list of the 29 step functions with five terms of Vasyunin. We also reformulate and generalize a conjecture of Vasyunin.
\end{abstract}
\email{borisov@math.psu.edu}
\maketitle
\section{Introduction}

The paper is organized as follows. In section 2 we discuss some known results on the integer ratios of factorials, mostly due to Landau and Picon. In section 3 we discuss the known classification results on cyclic quotient singularities and their connection to the integer ratios of factorials. In section 4 we discuss the Nyman's approach to the Riemann Hypothesis. In particular, we discuss some results of Vasyunin and their relationship with the other two topics of the paper. Section 5 is devoted to open questions and speculations.

{\bf Disclaimer.} The author was drawn into this circle of ideas from his work on the quotient singularities. So the paper is certainly imbalanced in that direction. This preprint is a preliminary version, and all suggestions and references will be greatly appreciated.

{\bf Acknowledgments.} The author thanks George Andrews, Ira Gessel, Andrew Granville and Gregory Sankaran for some helpful discussions related to this paper. The author originally learned of the Nyman's theorem from R. C. Vaughan.

\section{Integer Ratios of Factorials}

We start by restating two classical facts about the integer ratios of factorials. The first one is truly well known, while the second one is a little less popular.

\begin{Theorem} For all $n, k \in \N$
$$(n+k)!/(n!k!)\in \N .$$
\end{Theorem}

The integers above are usually called the binomial coefficients and have a very simple combinatorial interpretation.

\begin{Theorem} For all $n, k \in \N$
$$\frac{(2n)!(2k)!}{n!k!(n+k)!} \in \N .$$
\end{Theorem}

This theorem goes back at least to Catalan. However, no good combinatorial interpretation of these numbers is known, except when $k\leq 2$ (according to Ira Gessel).

The easiest way to prove this theorem is  by looking at the power with which an arbitrary prime number $p$ divides this ratio of factorials. This power equals the following formally infinite sum.
$$\sum \limits_{i=1}^{\infty} ([\frac{2n}{p^i}]+[\frac{2k}{p^i}]-[\frac{n}{p^i}]-[\frac{k}{p^i}]-[\frac{n+k}{p^i}])$$
One can check that for all real numbers $\alpha$ and $\beta$ 
$$[2\alpha]+[2\beta]-[\alpha]-[\beta]-[\alpha+\beta]\geq 0 \hskip 4cm (*)$$
Applying this to $\alpha=n/p^i$ and $\beta=k/p^i$ proves the result.

The main general result in this direction was published by Landau in 1900.

\begin{Theorem} (Landau, \cite{Landau}) Let $S$ be a finite set of linear forms on  ${\R}^d$ with coefficients in $\N$. To every $u\in S$ associate a unit $\varepsilon_u=\pm 1.$ For $X=(x_1,...,x_d) \in {\R}^d$ set $A(X)=\prod_{u\in S} ([u(X)]!)^{\varepsilon_u}$ and set $B(X)=\sum_{u\in S} [u(X)]{\varepsilon_u}.$ Then $A(X)\in{\N}$ for all $X\in {\N}^d$ iff $B(X)\geq 0$ for all $X \in [0,1]^d.$
\end{Theorem}

Using the above theorem, Landau proved, in particular, the following.

\begin{Theorem} For all $n, k \in \N$
$$\frac{(4n)!(4k)!}{(2n+k)!(n+2k)!n!k!} \in \N .$$
\end{Theorem}

This was developed further by P. A. Picon, who proposed several ways to construct new integer ratios of factorials out of the old ones (cf. \cite{Picon1}, \cite{Picon2}, \cite{Picon3}). In particular, he proved the following two theorems.

\begin{Theorem} (P.A. Picon) For all $n, k \in \N$
$$\frac{(2n)!(4k)!(4n+2k)!}{(n!)^2(k!)^2(n+k)!^2(2n+2k)!}\in \N .$$
\end{Theorem}
 
\begin{Theorem} (P.A. Picon) For all $n, k \in \N$
$$\frac{(12n)! n!(12k)!k!}{(3n!)(4n!)(3k)!(4k)!(n+k)!^6} \in \N .$$
\end{Theorem}

\centerline{\bf Cyclic Quotient Singularities}

A cyclic quotient singularity of dimension $d$ is an affine variety obtained by taking a quotient of the affine space ${\C}^d$ by a linear action of the group of $\{\mu_l\}$ of the $l-$th roots of unity. Diagonalizing, we can assume that this action is given as follows. For any $\mu \in \{\mu_l\}$ and $(x_1,x_2,...,x_d)\in {\C}^d$ we have $\mu (x_1,x_2,...,x_d)=(x_1 \mu^{a_1},x_2 \mu^{a_2},...,x_d \mu^{a_d})$ where $a_i$ are some integers from $0$ to $(l-1).$ In what follows we will assume that they are all non-zero, so that the singularity is not a Cartesian product of ${\C}$ and a lower-dimensional singularity. We will also assume that $l$ is the least common multiple of all $a_i$ so that the action is faithful. We will further assume that $l$ is the least common multiple for any collection of $a_i$ that misses one of the numbers. This means that the quotient map does not factor through the map obtained by raising one of the coordinates to some power. One can associate to this singularity a point $u =(a_1/l,a_2/l, ...a_d/l)$ in the ``standard torus" $T^d=\{(u_1,u_2,...,u_d) \in {\R}^d | 0\leq u_i <1\}$. Conversely, any such point defines a cyclic quotient singularity.

One of the most important invariants of a cyclic quotient singularity (for the purposes of the birational geometry) is its  minimal log-discrepancy. There are actually two versions of it, as the following definition explains.

\begin{Definition} Suppose a point $u =(a_1/l,a_2/l, ..., a_d/l)$ in  $T^d$ defines a cyclic quotient singularity of order $l$ of dimension $d.$ For $m=1,2,...,l-1$ we call the points $(\{m\cdot a_1/l\},\{m\cdot a_2/l\}, ..., \{m\cdot a_d/l\})$ the multiples of $u$. The Shokurov minimal log-discrepancy of the corresponding singularity is the minimal sum of the coordinates of those multiples that lie strictly inside $T^d.$ The total minimal log-discrepancy is the smallest sum of the coordinates of all points in the lattice generated by $u$ and ${\Z}^d$ with all coordinates being non-negative and at least two coordinates being positive.
\end{Definition}

The above definitions for the minimal log-discrepancy are the combinatorial interpretations for the cyclic quotient singularities of the general invariants that can be defined for all singularities (cf. \cite{Reid}, \cite{FlipsAbundance}). One can see immediately that the total minimal log-discrepancy is at most $2$, while the Shokurov minimal log-discrepancy is at most $d.$ The Shokurov minimal log-discrepancy is always greater than or equal to the total one. So it is a more subtle invariant of a singularity.
We are not going to use the total minimal log-discrepancy much in the paper, so for us the minimal log-discrepancy will always mean the Shokurov minimal log-discrepancy, unless explicitly specified otherwise.

The minimal log-discrepancy is always positive (for the quotient singularities). The singularity is called terminal (resp. canonical) if its (Shokurov {\it or} total) minimal log-discrepancy is greater than (resp. greater than or equal to) $1.$ This is equivalent to all non-zero multiples of the point $u=(a_1/l,a_2/l, ...a_d/l)$ in $T^d$ being outside of the simplex $\sum _i (u_i) \leq 1$ (resp. $<1$).

We now give a short overview of the known results on the classification of the terminal and canonical cyclic quotient singularities. More information on this and on the related topics can be found in the short survey \cite{ToricSurvey}.

Terminal and canonical singularities are the simplest kinds of singularities that have to be allowed in the higher dimension in order to generalize the classical birational geometry of smooth surfaces (cf., e.g., \cite{CKM}). In dimension two, there really are no terminal singularities (i.e. all of them are just  smooth), and the canonical singularities are the DuVal singularities. The canonical cyclic quotient singularities in dimension $2$ are the DuVal singularities of type $A$. They occur when the corresponding point in the torus $T_2$ lies on the segment $u_1+u_2=1.$ This segment is a subtorus of $T_2$ which lies outside of the interior of the standard simplex. 

If the dimension $d$ is $3$, then the classification of terminal cyclic quotient singularities is well known. Up to the permutation of the coordinates, the corresponding point $(a_1/l,a_2/l,a_3/l)$ lies in the subtorus of $T_3$ given by the hyperplane $u_1+u_2=1$. This is directly related to the classification of lattice tetrahedra with no lattice points inside, which was done by White back in 1952 (\cite{White}). A much shorter proof of it, using Bernoulli functions, was  given by D. Morrison and Stevens (\cite{MorStev}). The classification of the canonical cyclic quotient singularities in dimension 3 is also known (cf. \cite{Morrison}, \cite{IshidaIwashita}). Almost all of them are either terminal or correspond to the points $u$ with $u_1+u_2+u_3=1.$  

If $d=4$, the situation is more complicated. In 1988 Mori, Morrison and Morrison numerically classified all 4-dimensional terminal quotient singularities, with $l$ for simplicity being prime (\cite{MMM}). Instead of working with four rational numbers $a_1/l, a_2/l, a_3/l, a_4/l$ it turned out beneficial to work with five rational numbers $a_1/l, a_2/l, a_3/l, a_4/l, -(a_1+a_2+a_3+a_4)/l$ that sum up to zero. More precisely, for every 5-tuple of integers $(a_1,a_2,a_3,a_4,a_5)$ that sum up to zero and a (prime) number $l$ one can associate five cyclic quotient singularities by dropping one of the numbers and considering the point in $T^4$ with coordinates $\{a_i/l\}$ for the remaining four numbers. Each of these five singularities is terminal if and only if the quintuple satisfies the following condition.

{\bf Condition.} (Terminality condition) For all integers $k$ from $1$ to $l-1$
$$\sum \limits_{i=1}^{5} \{\frac{a_i}{l}\cdot k\} \geq 2 .$$
(Notice that this sum is always an integer.)

The main empirical result of the Mori-Morrison-Morrison paper is that apart from the large but finite number of the ``sporadic" terminal singularities for $l<421,$ all terminal cyclic quotient singularities belong to one of the three cases below.

1) The quintuple is of the form $(\alpha, -\alpha, \beta, \gamma, -\beta -\gamma)$

2) The quintuple is of the form $(\alpha, -2\alpha, \beta, -2\beta, \alpha+\beta)$

3) The quintuple belongs to one of the 29 ``exceptional stable quintuples" below. These are the quintuples of integers, two of which are positive and three are negative, (or the other way around) that work for any large enough $l$. For example, $(9,1,-2,-3,-5)$ will produce five terminal cyclic quotient singularities for every $l\geq 11$:
$$(\frac{9}{l}, \frac{1}{l}, \frac{l-2}{l}, \frac{l-3}{l});\ \   (\frac{1}{l}, \frac{l-2}{l}, \frac{l-3}{l}, \frac{l-5}{l});  \ \ (\frac{l-2}{l}, \frac{l-3}{l}, \frac{l-5}{l}, \frac{9}{l});\hskip 2cm$$
$$(\frac{l-3}{l}, \frac{l-5}{l}, \frac{9}{l}, \frac{1}{l}); \ \ (\frac{l-5}{l}, \frac{9}{l}, \frac{1}{l}, \frac{l-2}{l})\hskip 10cm.$$

The list of the 29 exceptional stable quintuples deserves to be reproduced here. We changed the notation of \cite{MMM} in the following way. For those quintuples that originally had three positive numbers, all numbers were multiplied by $-1.$ Then the two positive numbers are written first, separated by the semicolon from the absolute values of the three negative numbers.

\centerline{\bf Exceptional Stable Quintuples of Mori-Morrison-Morrison}

$$(9, 1; 2, 3, 5 ), \ (9, 2; 1, 4, 6), \ (12, 4; 4, 5,6), \ (12, 2; 3, 4, 7), \ (9, 4; 2, 3, 8),\hskip 5cm$$

$$(12, 1; 2, 3, 8), \ (12, 3; 1, 6, 8), \ (15, 4; 5, 6, 8), \ (12, 2; 1, 4, 9), \ (10, 6; 2, 5, 9),$$

$$(15, 1; 2, 5, 9), \ (12, 5; 3, 4, 10), \ (15, 2; 3, 4, 10), \ (12, 1; 3, 4, 6), \ (14, 1; 3, 5, 7),$$

$$(14, 3; 1, 7, 9), \ (15, 7; 3, 5, 14), \ (15, 1; 3, 5, 8), \ (15, 2; 1, 6, 10), \ \ (15, 4; 2, 5, 12),$$

$$(18, 1; 4, 6, 9), \ (18, 2; 5, 6, 9), \ (18, 4; 1, 9, 12), \ (20, 1; 4, 7, 10), \ (20, 1; 3, 8, 10),$$

$$(20, 3; 4, 9,10), \ (20, 3; 1, 10, 12), \ (24, 1; 5, 8, 12), \ (30, 1; 6, 10, 15). \hskip 10cm$$

It was proven by Gregory Sankaran \cite{Sankaran} that this list of the 29
stable quintuples is complete.

\vskip 1cm

The most general qualitative  result in this direction was obtained in \cite{vasq}. It implies in particular that in any fixed dimension the cyclic quotient singularities with Shokurov minimal log-discrepancy greater than (greater than or equal to) a fixed number always form a finite number of families. The main ingredient of the proof is the following beautiful theorem of Jim Lawrence (\cite{Lawrence}):

\begin{Theorem} (J. Lawrence, \cite{Lawrence})  Suppose $T$ is a torus (or, more generally, a closed subgroup of some torus). Suppose $U$ is an open subset of $T$ (or, more generally, a {\it full} subset of $T$, i.e. for all closed subgroups $L$ of $T$ the intersection of $L$ and $U$ is either empty or contains a relatively open subset of $L$). Consider all closed subgroups of $T$ that don't intersect $U$. Then the set of maximal elements of it, with respect to the inclusion, is finite.
\end{Theorem}

Together with the result of Sankaran, this implies that the numerical classification of  Mori-Morrison-Morrison is complete, except possibly for the finite number of the ``sporadic" singularities.

\vskip 1cm

Besides the minimal log-discrepancy, another important invariant of a  singularity is its Gorenstein index. In the case of the cyclic quotient singularity it equals to the least common multiple of the denominators of the sums of all coordinates of the corresponding points in $T^d$. One can further conclude that if the singularity is given by the point $(a_1,a_2,...,a_d)\in T^d$ then the Gorenstein index is the just the denominator of $\sum \limits_{i=1}^{n} a_i.$ The following class of singularities is very important.

\begin{Definition} The singularity is called Gorenstein if its Gorenstein index is $1.$
\end{Definition}

It is essentially obvious that the denominator of the minimal log-discrepancy divides the Gorenstein index of a singularity. Therefore all Gorenstein cyclic quotients have minimal log-discrepancy at least $1,$ and thus are canonical. 
(In fact, any log-terminal Gorenstein singularity is canonical).

The following observation shows that the Gorenstein cyclic quotient singularities occupy  a very special place among all cyclic quotients. This theorem is certainly not new, but we do not know if it was explicitly stated before.

\begin{Theorem} Suppose  the point $(a_1,a_2,...,a_d)\in T^d$ defines a cyclic quotient singularity, which is {\bf not} Gorenstein. Then the point 
$$(a_1,a_2,...,a_d, \{-\sum a_i\})\in T^{d+1}$$
defines a Gorenstein cyclic quotient singularity of dimension one bigger.  In the other direction, dropping a coordinate from a Gorenstein cyclic quotient produces a non-Gorenstein cyclic quotient of dimension one less. The minimal discrepancy of a non-Gorenstein cyclic quotient is less than the minimal discrepancy of the corresponding Gorenstein one by a rational number, which is less than one.
\end{Theorem}

{\bf Proof.} This follows from the inequality $0\leq \{-m\sum a_i\}<1.$

While not explicit in \cite{vasq}, the Gorenstein cyclic quotients of fixed dimension and minimal log-discrepancy also form a finite number of families. In view of the above theorem, it is most natural to interpret the results of Mori-Morrison-Morrison as the numerical classification of the $5-$dimensional Gorenstein cyclic quotients of prime index with minimal log-discrepancy $2$.

The following theorem connects the families of Gorenstein cyclic quotient singularities with the integer ratios of factorials.

\begin{Theorem} Suppose $u_1,u_2,...,u_k$ and $v_1,v_2,...,v_m$ are two finite sets of linear forms on  ${\R}^d$ with coefficients in $\N$. Suppose further that $\sum \limits_{i=1}^{k} u_i(X) = \sum \limits_{i=1}^{m} v_i(X)$ for all $X=(x_1,...,x_d) \in {\R}^d.$  Then the following two statements are equivalent.

\begin{enumerate}
\item For every $X=(x_1,...,x_d) \in {\N}^d$ 
$$\frac{\prod \limits_{i=1}^{k} u_i(X)!}{\prod \limits_{i=1}^{m} v_i(X)!} \in {\N}$$

\item For any $l\in {\N}$ and all $X=(x_1,...,x_d) \in {\Z}^d$ such that all $\{\frac{u_i(X)}{l}\}$ and $\{\frac{v_i(X)}{l}\}$ are non-zero, the following point in $T^{n+k}$ defines a Gorenstein cyclic quotient with Shokurov minimal log-discrepancy at least $k.$
$$\Big( \{\frac{-u_1(X)}{l}\},\{\frac{-u_2(X)}{l}\},...,\{\frac{-u_k(X)}{l}\},\{\frac{v_1(X)}{l}\},\{\frac{v_2(X)}{l}\},...,\{\frac{v_m(X)}{l}\}\Big) $$ 
\end{enumerate}

\end{Theorem}

{\bf Proof.} First we will prove that $(1)$ implies $(2)$. By the theorem of Landau (Theorem 3) $(1)$ implies that the following inequality is satisfied for all $l \in {\N}$ and all $X$ in ${\R}^d.$

$$\sum \limits_{i=1}^{k} [u_i(\frac{X}{l})] \geq \sum \limits_{i=1}^{m} [v_i(\frac{X}{l})]$$

Because $\sum u_i(\frac{X}{l})= \sum v_i(\frac{X}{l}),$ this implies that

$$ \sum \limits_{i=1}^{k} \{u_i(\frac{X}{l})\} \leq \sum \limits_{i=1}^{m} \{v_i(\frac{X}{l})\}.$$

If all $\{ u_i(\frac{-X}{l}) \} \neq 0,$ then we can use the identity $\{-x\}=1-\{x\}$ to get
$$k- \sum \limits_{i=1}^{k} \{u_i(\frac{-X}{l})\} \leq \sum \limits_{i=1}^{m} \{v_i(\frac{X}{l})\},$$ which is equivalent to

$$\sum \limits_{i=1}^{k} \{u_i(\frac{-X}{l})\} +\sum \limits_{i=1}^{m} \{v_i(\frac{X}{l})\} \geq k.$$

The result now follows from the definition of the Shokurov minimal log-discrepancy (see Definition 7).

Now we prove that $(2)$ implies $(1)$. This requires a more delicate argument because in the definition of the Shokurov minimal log-discrepancy we only use the multiples that lie directly inside $T.$ It is conceivable that some multiples on the boundary have a smaller sum of the coordinates.

We denote by $\bar{\{} x \bar{\}}$ the ``upper fractional part of $x$". It equals $\{x\}$ if $x$ is not an integer and $1$ if it is. For all $x$ we have $\bar{\{} x \bar{\}}=1-\{-x\}.$

Define the function $F(X)$ on $X=(x_1,...,x_d)\in {\R}^d$ as follows:
$$F(X)=\sum \limits_{i=1}^{k} \bar{\{}u_i(-X)\bar{\}} +\sum \limits_{i=1}^{m} \{v_i(X)\}$$

Clearly, for the rational points $X$ such that all $u_i(X)$ and $v_i(X)$ are not integers, $F(X)$ is greater than or equal to the minimal log-discrepancy of the singularity that corresponds to the point
$$ \{-u_1(X)\},\{-u_2(X)\},...,\{-u_k(X)\},\{v_1(X)\},\{v_2(X)\},...,\{v_m(X)\}.$$

Thus for such points $X$ we have $F(X)\geq k$. One can also note that for all $X$ 
$$F(X)=\lim \limits_{\varepsilon _i \to 0+} F(x_1+\varepsilon _1, x_2+\varepsilon _2,...,x_d+\varepsilon _d) $$ 
This is because the coefficients of $u_i$ and $v_i$ are positive and $\{x\}$ and $\bar{\{}x\bar{\}}$ are right-continuous and left-continuous respectively.

The above limit implies that $F(X)\geq k$ for all $X\in {\R}^d.$ Because ${u}=1-\bar{\{}-u \bar{\}}$ for all $u,$ this implies that for all $X\in {\R}^d$

$$ \sum \limits_{i=1}^{k} \{u_i(X)\} \leq \sum \limits_{i=1}^{m} \{v_i(X)\}.$$

Because $\sum u_i(\frac{X}{l})= \sum v_i(\frac{X}{l}),$ this implies that

$$ \sum \limits_{i=1}^{k} [u_i(X)] \geq \sum \limits_{i=1}^{m} [v_i(X)].$$

The result now follows from the Landau's Theorem (Theorem 3).

\vskip 1cm

The above theorem allows us to translate the results of Mori-Morrison-Morrison into the statements on the integrality of some ratios of factorials. For example, the inequalities in the second case above lead to the numbers of Catalan. Each of the 29 stable quintuples in the third case leads to a family of integer ratios of factorials with one parameter. For example, the quintuple $(9,1;2,3,5)$ leads to the sequence
$$\frac{(9n)!n!}{(5n)!(3n)!(2n)!}.$$ 
We also have 
$$\frac{(14n)!(3n)!}{(9n)!(7n)!n!}, \frac{(14n)!n!}{(7n)!(5n)!(3n)!}, \frac{(20n)!(3n)!}{(10n)!(9n)!(4n)!}$$
and $25$ more examples of the same nature.

The result of Sankaran implies that every family of the integer ratios of the factorials of the form $\frac{(an)!(bn)!}{(cn)!(dn)!(en)!}$ comes from one of these 29 families or from the two-parameter family $\frac{(2n)!(2k)!}{(n+k)!n!k!}$ by means of a linear substitution. 

In the other direction, the results of Picon provide some new and interesting cyclic quotient singularities.

\section{Nyman's approach to the Riemann Hypothesis}
The binomial coefficients $\frac{(2n)!}{(n)!^2}$ have been used in the elementary estimates for the number of primes (Chebyshev). So it is not entirely surprising that the more intricate ratios of factorials can be used for the more detailed study of the distribution of primes. The principal result in this direction is the following celebrated theorem of Nyman (\cite{Nyman}).

\begin{Theorem} (Nyman-Beurling Criterion for RH) The Riemann Hypothesis is equivalent to the following statement.
The set of linear combinations of the following functions is dense in $L^2([0,1]):$
$$\rho_{\alpha}(x)=[\frac{\alpha}{x}]-\alpha[\frac{1}{x}], \ \ 0<\alpha <1.$$
\end{Theorem}

In fact, the statement above can be replaced by the possibility of approximating by the linear combinations of $\rho_{\alpha}$ a single function: the constant $1.$  Several further improvements were obtained. In particular,  B\'aez-Duarte proved (cf. \cite{Baez-Duarte}) that in the above statement one can restrict $\alpha$ to be of the form $1/n, \ n\in {\N}.$ 

It is actually convenient to substitute $\frac{1}{x}$ for $x$ in the Nyman-Beurling criterion and consider the functions on $[1,\infty)$ instead of $[0,1].$ The relevant $L^2$ space will be with respect to the measure $\frac{1}{x^2} dx.$ The functions $\rho$ for $\alpha=1/n$ will look like this:

$$\rho_{n}(x)=[\frac{x}{n}]-\frac{1}{n}[x], \ \ n\in {\N}.$$
These functions are piecewise constant, with jumps at the natural numbers that are multiples of $n.$ Their finite linear combinations can be described as the following sums over the finite collection $S$ of natural numbers $n$.
$$F(S,c_n;x)= \sum \limits_{n\in S} c_n [\frac{x}{n}] \ {\textrm where} \  \sum \limits_{n\in S} \frac{c_n}{n} =0 .$$

It is obviously sufficient to require the coefficients $c_n$ to be rational numbers. Multiplying $F(S,c_n;x)$ by the common denominator of $c_n$ and approximating some bigger constant instead of the constant $1,$ one gets the linear combinations of $[\frac{x}{n}]$ with integer coefficients. As explained in the previous two sections, they are directly related to the integer ratios of factorials and the cyclic quotient singularities. One has to note that the singularities considered are of increasingly large dimension. 

Of particular interest to us is the paper of Vasyunin \cite{Vasyunin}. He was motivated by the idea of modifying the functions $F(S,c_n;x)$ to achieve a good approximation for small $x$ while trying not to spoil the behavior for large $x$ too much. This led him to the problem of classifying all functions   $F(S,c_n;x)$ that take only values $0$ and $1.$ His first result is the following.

\begin{Theorem}(Vasyunin, \cite{Vasyunin}) Suppose  $F(S,c_n;x)$ is as above. Suppose it only takes values $0$ and $1$ for $x\geq 1$. Then we have the following
\begin{enumerate}

\item All coefficients $c_n$ are integers.

\item $\sum _n c_n=-1.$ This means that we can rewrite $F(S,c_n;x)$ as $$\sum \limits_{i=1}^{d} [\frac{x}{n_i}]- \sum \limits_{i=d+1}^{2d+1} [\frac{x}{n_i}] $$ with not necessarily distinct $n_i.$
\end{enumerate}

\end{Theorem}

In the formula above one can require that no numbers $n_i$ for $i=1,2,..,d$ equal some numbers $n_i$ for $i=d+1,...,2d+1$. Vasyunin calls such expressions the irreducible solutions to his problem. Clearly, these are the only solutions worthy of consideration. Additionally, one can always multiply or divide $n_i$ by a common factor, which corresponds to the scaling of the variable $x.$

Together with the discussion in the previous section, the theorem above shows that classifying the functions  $F(S,c_n;x)$ with the desired properties is equivalent to classifying the one-dimensional families of Gorenstein quotient singularities of dimension $2d+1$ with Shokurov minimal log-discrepancy $d.$ Note that because of the natural symmetry this is the biggest minimal log-discrepancy possible for the $(2d+1)-$dimensional cyclic quotient singularities.

For $d=1$ all Gorenstein quotient singularities of dimension $3$ have minimal log-discrepancy $1.$ So we get a function $F(S,c_n;x)$ with the desired property for all triples of integers $(n_1;n_2,n_3)$ such that $\frac{1}{n_1}=\frac{1}{n_2}+\frac{1}{n_3}.$  Vasyunin was unaware of the geometric interpretation of his problem, so for him this classification is just an elementary exercise.

For $d=2$ the problem of Vasyunin is equivalent, as explained in the previous section, to the problem considered by Mori-Morrison-Morrison in \cite{MMM} (disregarding the sporadic singularities). Vasyunin also used computers to get numerically the classification of the ``five term step-functions"  $F(S,c_n;x)$ that only take values $0$ and $1$. Sure enough, his results are in perfect agreement with those of Mori-Morrison-Morrison. The notation is different, but one can easily match his 29 special irreducible solutions with the 29 series of Mori-Morrison-Morrison. Vasyunin conjectures that this is a complete list. He actually sketches the proof of this, but notes that it is unacceptably tedious. The proof, which takes a considerable part of the paper \cite{Vasyunin}, seems to be close in the spirit to the proof of Sankaran, of which Vasyunin is obviously unaware. The proof of Sankaran is actually so long that less than half of it was formally published (cf. \cite{Sankaran}). The argument of Vasyunin seems shorter but harder to check. It is not clear whether deep down the two proofs are essentially the same.

Vasyunin's computer experiments went beyond those of \cite{MMM}. For $d=3$ he found 21 functions $F(S,c_n;x)$, and for $d=4$ he found 2. Thus he found 21 one-parameter families of $7-$dimensional cyclic quotients with Shokurov minimal log-discrepancy $3$ and 2 families in dimension $9$ with minimal log-discrepancy $4.$ It seems like they were not considered before.

Vasyunin could not prove the classification for $d>2$ so he just conjectured that his lists for $d=3$ and $d=4$ are complete. Most remarkably, his computer experiments show that there are no functions  $F(S,c_n;x)$ that only take values $0$ and $1$ for $d\geq 5.$ He conjectures this to be true and speculates that there should be some simple arithmetical reason behind this.

\section{Conjectures and Speculations}

In the language of cyclic quotient singularities, the conjecture of Vasyunin is equivalent to the following.

\begin{Conjecture} Suppose $d\geq 5$ and we have a one-parameter family of Gorenstein cyclic quotient singularities of dimension $2d+1$ with Shokurov minimal log-discrepancy $d.$ Then up to the permutation of the coordinates in the $T^{(2d+1)},$ the corresponding points lie in the subtorus $x_1+x_2=1.$
\end{Conjecture}

One can see that the corresponding result for the one-dimensional families of  the Gorenstein singularities of dimension $2d$ and minimal log-discrepancy $d$ is equivalent to the numbers $n_i$ being uniquely determined by the function $F(S,c_n;x) \sum \limits_{n\in S} c_n [\frac{x}{n}].$ This is easy to check by looking at the jumps of the step-function. For the individual singularities this is only proven,  to the author's knowledge, for $d=1$ (trivial) and $d=2$ (the Terminal Lemma of White and Morrison-Stevens).

Based more on our intuition than on any hard evidence, we hereby boldly propose the following far-fetching strengthening of the conjecture of Vasyunin.

\begin{Conjecture} (Higher-dimensional Terminal Lemma)
Suppose $a\geq 0$ is any real number. Then for all big enough $d\in {\N}$ the following statement is true.

If $(x_1,...,x_d)\in T^d$ defines a cyclic quotient singularity with Shokurov minimal log-discrepancy at least $\frac{d}{2}-a,$ then for some pair of indices $1\leq i <j\leq d$ we have $x_i+x_j=1.$
\end{Conjecture}

We do not know at the moment if this conjecture has any implications for the distribution of primes, but it is conceivable that it does. The most promising approach to proving it would be to generalize the Bernoulli function argument of \cite{MorStev}.

In \cite{ToricSurvey} we conjectured that one can automatize the argument of Sankaran. No one has done it yet, due to the lack of a serious incentive. Hopefully this paper will enhance the interest in this topic.

Another interesting possibility is to try to extend the result of Lawrence \cite{Lawrence} to the situation when the subgroups intersect with the open set, but the intersection is in some sense small. One would also like to investigate whether any of the known discrete geometry results (cf. \cite{ToricSurvey}) have any implications for the distribution of primes.

Ultimately, to advance the knowledge of the distribution of primes via the Nyman's approach one should learn to construct many ``good" step-functions (or, in other words, ``good" cyclic quotient singularities. A few such constructions already exist (cf. \cite{Toricdiscr}, \cite{Picon1}) but much more can be done in this direction.

\end{document}